\pdfoutput=1
\RequirePackage{ifpdf}
\ifpdf
\documentclass[pdftex]{sigma}
\else
\documentclass{sigma}
\fi

\usepackage{enumerate}

\DeclareMathOperator{\Cliff}{Clif\/f}
\DeclareMathOperator{\Gal}{Gal}
\DeclareMathOperator{\GL}{GL}

\DeclareMathOperator{\Mat}{Mat}

\DeclareMathOperator{\SU}{\mathbf{SU}}
\numberwithin{equation}{section}

\newtheorem{prop}{Proposition}[section]
\newtheorem{cor}[prop]{Corollary}
\newtheorem{lem}[prop]{Lemma}
{\theoremstyle{definition}
\newtheorem{defn}[prop]{Definition}
\newtheorem{notation}[prop]{Notation}
\newtheorem{rem}[prop]{Remark}

\newtheorem{eg}[prop]{Example}
\newtheorem{erratum}[prop]{Erratum}}

\begin{document}

\newcommand{\Res}[1]{\mathop{\mathbf{R}_{#1}}}

\newcommand{\ResOne}[1]{\mathop{\mathbf{R}^{(1)}_{#1}}}

\newcommand{\Lie}[1]{\mathfrak{\lowercase{#1}}}

\newcommand{\arXivNumber}{1410.2339}

\allowdisplaybreaks

\renewcommand{\PaperNumber}{034}

\FirstPageHeading

\ShortArticleName{Real Representations of Semisimple Lie Algebras Have $\mathbb{Q}$-Forms}

\ArticleName{A Cohomological Proof that Real Representations\\
of Semisimple Lie Algebras Have $\boldsymbol{\mathbb{Q}}$-Forms}

\Author{Dave Witte MORRIS}

\AuthorNameForHeading{D.W.~Morris}

\Address{Department of Mathematics and Computer Science, University of Lethbridge,\\
Lethbridge, Alberta, T1K~3M4, Canada}
\Email{\href{mailto:Dave.Morris@uleth.ca}{Dave.Morris@uleth.ca}}
\URLaddress{\url{http://people.uleth.ca/~dave.morris/}}

\ArticleDates{Received October 17, 2014, in f\/inal form April 14, 2015; Published online April 27, 2015}

\Abstract{A Lie algebra $\Lie G_\mathbb{Q}$ over~$\mathbb{Q}$ is said to be \emph{$\mathbb{R}$-universal} if every
homomorphism from~$\Lie G_\mathbb{Q}$ to $\Lie{GL}(n,\mathbb{R})$ is conjugate to a~homomorphism into
$\Lie{GL}(n,\mathbb{Q})$ (for every~$n$).
By using Galois cohomology, we provide a~short proof of the known fact that every real semisimple Lie algebra has an
$\mathbb{R}$-universal $\mathbb{Q}$-form.
We also provide a~classif\/ication of the $\mathbb{R}$-universal Lie algebras that are semisimple.}


\Classification{17B10; 17B20; 11E72; 20G30}

\section{Introduction}

\begin{defn}
Let $\Lie G_\mathbb{Q}$ be a~Lie algebra over~$\mathbb{Q}$.
(All Lie algebras and all representations are assumed to be f\/inite-dimensional.)
\begin{enumerate}\itemsep=0pt
\item $\Lie G_\mathbb{Q}$ is \emph{universal for real representations} (or \emph{$\mathbb{R}$-universal}, for short) if
every real representation of~$\Lie G_\mathbb{Q}$ has a~$\mathbb{Q}$-form~\cite[Def\/inition~7.1]{Morris-RRepsQForms}.
This means that if $\rho \colon \Lie G_\mathbb{Q} \to \Lie{GL}(n,\mathbb{R})$ is any ($\mathbb{Q}$-linear) Lie algebra
homomorphism, then there exists $M \in \GL(n,\mathbb{R})$, such that $M \rho(x) M^{-1} \in \Lie{GL}(n,\mathbb{Q})$, for
every $x \in \Lie G_\mathbb{Q}$.
\item $\Lie G_\mathbb{Q}$ is a~\emph{$\mathbb{Q}$-form} of a~real Lie algebra~$\Lie G_\mathbb{R}$ if $\Lie G_\mathbb{Q}
\otimes_\mathbb{Q} \mathbb{R} \cong \Lie G_\mathbb{R}$.
\end{enumerate}
\end{defn}

This note uses Galois cohomology to present a~short proof of the following known result, which was f\/irst obtained~by
M.S.~Raghunathan~\cite[\S~3]{Raghunathan-ArithLattsSSGrps} in the important special case where $\Lie G_\mathbb{R}$ is
compact.

\begin{prop}[\protect{\cite[Theorem~1.2]{Morris-RRepsQForms}}]\label{HasRUniv}
Every real semisimple Lie algebra has a~$\mathbb{Q}$-form that is $\mathbb{R}$-universal.
\end{prop}

The proof in~\cite{Morris-RRepsQForms} constructs an $\mathbb{R}$-universal $\mathbb{Q}$-form explicitly, and is rather
tedious, but a~much nicer proof was given by G.~Prasad and A.~Rapinchuk~\cite[Proposition~3 and Remark~3]{PrasadRapinchuk-Prescribed}.
Assuming some fundamental results of J.~Tits~\cite{Tits-RepLinIrred}, our proof in Section~\ref{PfSect} is a~bit shorter and
more direct.
(On the other hand, we provide less information about the $\mathbb{Q}$-form than is supplied
in~\cite{PrasadRapinchuk-Prescribed}.)

Section~\ref{WhichRUnivSect} gives an explicit characterization of the $\mathbb{R}$-universal Lie algebras that are absolutely
simple over~$\mathbb{Q}$, and, for completeness, Section~\ref{NotAbsSect} explains how to extend this to the class of all
semisimple Lie algebras.

Due to the well-known correspondence between $\mathbb{Q}$-forms and arithmetic subgroups, Proposition~\ref{HasRUniv} has the
following consequence in the theory of discrete subgroups:

\begin{cor}
\label{LattAppl}
Let~$G$ be a~connected, semisimple Lie group with finite center.
Then there is a~discrete subgroup~$\Gamma$ of~$G$, such that
\renewcommand{\labelenumi}{{\rm (\theenumi)}}
\begin{enumerate}\itemsep=0pt
\item $G/\Gamma$ has finite volume $($so~$\Gamma$ is a~``lattice'' in~$G)$, and
\item if $\rho \colon G \to \GL(n,\mathbb{R})$ is any finite-dimensional representation of~$G$, then $\rho(\Gamma)$ is
conjugate to a~subgroup of~$\GL(n,\mathbb{Z})$.
\end{enumerate}
\end{cor}

\section{Proof of the results stated in the introduction}\label{PfSect}

We begin by recalling a~result of J.~Tits that uses Galois cohomology to characterize the irreducible representations of
semisimple algebraic groups over f\/ields that are not algebraically closed.
However, we will state this work in the setting of Lie algebras, rather than algebraic groups.
(We deal only with semisimple Lie algebras and semisimple groups, and the f\/ields under consideration in this paper are
all of characteristic zero, so no dif\/f\/iculties arise in making this translation~\cite[Proposition~7.3.1(iii), p.~393]{DemazureGrothendieck-SGA3}.)

\begin{defn}[\protect{\cite[\S~4.2]{Tits-RepLinIrred}}] \label{BetaDefn}
Suppose $\Lie G$ is a~semisimple Lie algebra over a~subf\/ield~$F$ of~$\mathbb{C}$, and let $\mathbf{G}$ be the
corresponding simply connected, semisimple algebraic group over~$F$.
It is well known that there is a~(unique) quasi-split, simply connected algebraic group $\mathbf{G}^q$ over~$F$, and
a~$1$-cocycle $\xi \colon \Gal(\overline{F}/F) \to \overline{\mathbf{G}}{}^q$, where $\overline{\mathbf{G}}{}^q$ is the
adjoint group of~$\mathbf{G}^q$, such that $\mathbf{G}$ is~$F$-isomorphic to the Galois twist ${}^\xi \mathbf{G}^q$.
This cocycle determines a~cohomology class $[\xi] \in H^1(F; \overline{\mathbf{G}}{}^q)$.

Letting $Z(\mathbf{G}^q)$ be the center of $\mathbf{G}^q$, the short exact sequence $ e \to Z(\mathbf{G}^q) \to
\mathbf{G}^q \to \overline{\mathbf{G}}{}^q \to e $ yields a~corresponding long exact sequence of Galois cohomology sets,
including a~connecting map $\delta_* \colon H^1(F; \overline{\mathbf{G}}{}^q) \to H^2\bigl(F; Z(\mathbf{G}^q) \bigr)$.
Hence, we have a~cohomology class $\delta_* [\xi] \in H^2\bigl(F; Z(\mathbf{G}^q) \bigr)$.
(If we took a~bit more care to ensure that it is well-def\/ined, this would be the \emph{Tits class}
of~$\mathbf{G}$~\cite[p.~426]{KnusEtAl-BookOfInvols}.)

Now, f\/ix a~maximal~$F$-torus~$\mathbf{T}$ of~$\mathbf{G}^q$ that contains a~maximal~$F$-split torus, and
suppose~$\lambda$ is a~weight of~$\mathbf{T}$ that is invariant under the $*$-action of the Galois group $\Gal(\overline{F}/F)$.
Then the restriction of~$\lambda$ to $Z(\mathbf{G}^q)$ is a~$\Gal(\overline{F}/F)$-equivariant homomorphism from
$Z(\mathbf{G}^q)$ to the group~$\boldsymbol\mu$ of roots of unity in~$\mathbb{C}$, so it induces a~homomorphism
$\lambda_* \colon H^2\bigl(F; Z(\mathbf{G}^q) \bigr) \to H^2\bigl(F; \boldsymbol\mu \bigr)$.
Therefore, we may def\/ine
\begin{gather*}
\beta_{\Lie G,F}(\lambda) = \lambda_* \delta_* [\xi] \in H^2(F; \boldsymbol\mu).
\end{gather*}
\end{defn}

\begin{rem}
For a~semisimple group~$\mathbf{G}$ that is def\/ined over a~f\/ield~$F$, the def\/inition of the $*$-action of
$\Gal(\overline{F}/F)$ on the weights of a~maximal torus~$\mathbf{T}$ presupposes that $\mathbf{T}$ is def\/ined over~$F$
and contains a~maximal~$F$-split torus of~$\mathbf{G}$~\cite[p.~66]{PlatonovRapinchukBook}.
When $\mathbf{G}$ is def\/ined over~$\mathbb{Q}$, we will use the $*$-actions of both $\Gal(\overline{\mathbb{Q}}/\mathbb{Q})$ and $\Gal(\mathbb{C}/\mathbb{R})$, so we will assume that $\mathbf{T}$ is def\/ined
over~$\mathbb{Q}$ (so it is also def\/ined over~$\mathbb{R}$) and contains both a~maximal $\mathbb{Q}$-split torus and
a~maximal $\mathbb{R}$-split torus.
To see that such a~choice of~$\mathbf{T}$ is always possible, let $\mathbf{T}_1$ be any maximal $\mathbb{Q}$-split torus
of~$\mathbf{G}$.
Then there is a~maximal $\mathbb{Q}$-torus~$\mathbf{T}$ of the centralizer $\mathbf{C}_{\mathbf{G}}(\mathbf{T}_1)$ that
contains a~maximal $\mathbb{R}$-split torus~\cite[Corollary~2 of Proposition~7.8, p.~418]{PlatonovRapinchukBook}, and it is clear
that $\mathbf{T}$ has the desired properties.
\end{rem}

\begin{prop}[\protect{\cite[Theorem~7.2 and Lemma~7.4]{Tits-RepLinIrred}}]\label{HighWt}
Suppose $\Lie G$ is a~semisimple Lie algebra over a~subfield~$F$ of~$\mathbb{C}$, and~$\lambda$ is a~dominant weight.
Then:
\renewcommand{\labelenumi}{{\rm \theenumi.}}
\begin{enumerate}\itemsep=0pt
\item There is an irreducible representation ${}^F \rho_\lambda \colon \Lie G \to \Lie{GL}(n, F)$, for some~$n$, such
that \mbox{${}^F \rho_\lambda \otimes_F \mathbb{C}$} has an irreducible summand with highest weight~$\lambda$.
Furthermore, ${}^F \rho_\lambda$ is unique up to isomorphism.

\item ${}^F \rho_{\lambda_1} \cong {}^F \rho_{\lambda_2}$ if and only if~$\lambda_1$ and~$\lambda_2$ are in the same
orbit of the $*$-action of~$\Gal(\overline{F}/F)$.

\item
\label{HighWt-IrredIff}
${}^F \rho_\lambda \otimes_F \mathbb{C}$ is irreducible if and only if:
\begin{enumerate}[{\rm (a)}]\itemsep=0pt
\item
$\lambda$ is invariant under the $*$-action of~$\Gal(\overline{F}/F)$, and
\item
$\beta_{\Lie G,F}(\lambda)$ is the trivial element of $H^2(F; \boldsymbol\mu)$.
\end{enumerate}
\end{enumerate}
Furthermore, every~$F$-irreducible representation of~$\Lie G$ is isomorphic to ${}^F \rho_{\lambda}$, for some
dominant weight~$\lambda$.
\end{prop}

\begin{cor}
Suppose~$\Lie G$ is a~semisimple Lie algebra over~$\mathbb{Q}$, such that~$\Lie G$ splits over a~quadratic
extension of~$\mathbb{Q}$, and, for every dominant weight~$\lambda$ of~$\mathbf{T}$:
\begin{gather}
\begin{split}
& \text{if~$\lambda$ is invariant under the $*$-action of~$\Gal(\mathbb{C}/\mathbb{R})$, and $\beta_{\Lie G,\mathbb{R}}(\lambda) = 0$,}
\\
& \text{then~$\lambda$ is also invariant under the $*$-action of $\Gal(\overline{\mathbb{Q}}/\mathbb{Q})$, and
$\beta_{\Lie G,\mathbb{Q}}(\lambda) = 0$.}
\end{split}\label{SplitOverQuadIff-weight}
\end{gather}
Then $\Lie G$ is $\mathbb{R}$-universal.
\end{cor}

\begin{proof}
Since representations of $\Lie G$ are completely reducible, it suf\/f\/ices to show that every \emph{irreducible} real
representation of~$\Lie G$ has a~$\mathbb{Q}$-form.
Specif\/ically, we will show that ${}^\mathbb{Q} \rho_\lambda$ is a~$\mathbb{Q}$-form of~${}^\mathbb{R} \rho_\lambda$, for
every dominant weight~$\lambda$.

We begin by showing that if ${}^\mathbb{R} \rho_\lambda \otimes_\mathbb{R} \mathbb{C}$ is irreducible, then
${}^\mathbb{Q} \rho_\lambda \otimes_\mathbb{Q} \mathbb{C}$ is irreducible.
From Proposition~\ref{HighWt}(\ref{HighWt-IrredIff}) (with $F = \mathbb{R}$), we know that~$\lambda$ is invariant under the $*$-action
of $\Gal(\mathbb{C}/\mathbb{R})$, and that $\beta_{\Lie G,\mathbb{R}}(\lambda) = 0$.
By assumption, this implies that~$\lambda$ is invariant under the $*$-action of\/
$\Gal(\overline{\mathbb{Q}}/\mathbb{Q})$, and that $\beta_{\Lie G,\mathbb{Q}}(\lambda) = 0$.
Then, from Proposition~\ref{HighWt}(\ref{HighWt-IrredIff}) (with $F = \mathbb{Q}$), we conclude that ${}^\mathbb{Q} \rho_\lambda
\otimes_\mathbb{Q} \mathbb{C}$ is irreducible, as desired.

Since $\Lie G_\mathbb{Q}$ splits over a~quadratic extension, we know that ${}^\mathbb{Q} \rho_\lambda \otimes_\mathbb{Q}
\mathbb{C}$ is either irreducible or the direct sum of two irreducibles~\cite[Corollary~3.2(2)]{Morris-RRepsQForms}.
Therefore, ${}^\mathbb{Q} \rho_\lambda \otimes_\mathbb{Q} \mathbb{C}$ and ${}^\mathbb{R} \rho_\lambda \otimes_\mathbb{R}
\mathbb{C}$ have the same number of irreducible constituents.
(Namely, either they are both irreducible, or they are both the direct sum of $2$ irreducibles.) Since ${}^\mathbb{R}
\rho_\lambda$ is a~summand of ${}^\mathbb{Q} \rho_\lambda \otimes_\mathbb{Q} \mathbb{R}$, this implies that
${}^\mathbb{R} \rho_\lambda \cong {}^\mathbb{Q} \rho_\lambda \otimes_\mathbb{Q} \mathbb{R}$, so ${}^\mathbb{Q}
\rho_\lambda$ is a~$\mathbb{Q}$-form of ${}^\mathbb{R} \rho_\lambda$.
\end{proof}

We will also use the following variant of a~basic fact in the theory of Galois cohomology:

\begin{lem}[\protect{\cite[Theorem~5.1b, p.~77]{Kneser-LectGalCohoClassical}}]\label{KneserHasse}
If $\mathbf{G}$ is a~connected, semisimple algebraic group over~$\mathbb{Q}$, and~$L$ is any imaginary quadratic
extension of~$\mathbb{Q}$, then the natural restriction map $H^1 \bigl(L/F; \mathbf{G}(L) \bigr) \to H^1 (\mathbb{R};
\mathbf{G})$ is surjective.
\end{lem}

\begin{proof}(The author thanks A.~Rapinchuk for suggesting this argument.) Let~$\sigma$ be the nontrivial element of~$\Gal(\mathbb{C}/\mathbb{R})$.
It is well known that any cohomology class in $H^1(\mathbb{R}; \mathbf{G})$ is represented by an element~$t$ of
a~maximal torus $\mathbf{T}(\mathbb{C})$ of $\mathbf{G}(\mathbb{C})$, such that $t t^\sigma = 1$, and, since the variety
of maximal tori has weak approximation, that $\mathbf{T}$ may be chosen to be def\/ined over~$\mathbb{Q}$ (cf.\ proof
of~\cite[Proposition~6.17, p.~337]{PlatonovRapinchukBook}).

Let $\mathbf{R} = \mathrm{Res}_{L/\mathbb{Q}} \mathbf{T}$ be the torus obtained from~$\mathbf{T}$ by restriction of
scalars from~$L$ to~$\mathbb{Q}$.
Then there is a~natural isomorphism $\varphi \colon \mathbf{T}(\mathbb{C}) \to \mathbf{R}(\mathbb{R})$, such that
$\varphi \bigl(\mathbf{T}(L) \bigr) = \mathbf{R}(\mathbb{Q})$, and there is a~$\mathbb{Q}$-automorphism~$\tau$
of~$\mathbf{R}$, such that $\varphi(x^\sigma) = \varphi(x)^\tau$, for all $x \in \mathbf{T}(\mathbb{C})$.

Let $\mathbf{R}^{(1)} = \{ r \in \mathbf{R} \mid r r^\tau = 1 \}$.
We claim that this is a~subtorus of~$\mathbf{R}$.
Since $\mathbf{T}$ is def\/ined over~$\mathbb{Q}$, there is an~$L$-isomorphism $\psi \colon \mathbf{R} \to \mathbf{T}
\times \mathbf{T}$, such that $\psi \bigl(\mathbf{R}(\mathbb{R}) \bigr) = \{ (x, x^\sigma) \mid x \in
\mathbf{T}(\mathbb{C}) \}$.
Therefore, if we def\/ine $(x,y)^\pi = (y,x)$, then $\psi(r^\tau) = \psi(r)^\pi$ for all $r \in \mathbf{R}$.
So $\psi(\mathbf{R}^{(1)}) = \{ (x, x^{-1}) \mid x \in \mathbf{T} \}$ is a~subtorus of $\mathbf{T} \times \mathbf{T}$,
as claimed.

Now, $\mathbf{R}^{(1)}$ is a~torus that is def\/ined over~$\mathbb{Q}$ (since the automorphism~$\tau$ is def\/ined
over~$\mathbb{Q}$), and all $\mathbb{Q}$-tori have weak approximation at the inf\/inite place~\cite[Corollary~1 of Proposition~7.8,
p.~418]{PlatonovRapinchukBook}, so $\mathbf{R}^{(1)}(\mathbb{Q})$ is dense in $\mathbf{R}^{(1)}(\mathbb{R})$.
Hence, some $q \in \mathbf{R}^{(1)}(\mathbb{Q})$ is in the same connected component of $\mathbf{R}^{(1)}(\mathbb{R})$ as
$\varphi(t)$.
Then, letting $q' = \varphi^{-1}(q)$, we see that $q'$ is in the same connected component of $\{ w \in
\mathbf{T}(\mathbb{R}) \mid w w^\sigma = 1 \}$ as~$t$, so $q'$ represents the same cohomology class as~$t$ in
$H^1(\mathbb{R}; \mathbf{G})$.
However, since $q' \in \mathbf{T}(L)$, we see that the cohomology class of~$q'$ is in the image of $H^1 \bigl(L/F;
\mathbf{G}(L) \bigr)$, as desired.
\end{proof}

\begin{proof}[\bf Proof of Proposition~\ref{HasRUniv}] Suppose $\Lie G_\mathbb{R}$ is a~real semisimple Lie algebra, and let $\mathbf{G}$ be the
simply connected, semisimple $\mathbb{R}$-algebraic group whose Lie algebra is~$\Lie G_\mathbb{R}$.
As in Def\/inition~\ref{BetaDefn}, write $\mathbf{G} = {}^\xi \mathbf{G}^q$, where $\xi \colon \Gal(\mathbb{C}/\mathbb{R}) \to
\overline{\mathbf{G}}{}^q$ is a~$1$-cocycle and $\mathbf{G}^q$ is quasi-split.
Let $L = \mathbb{Q}[i]$.
By choosing an appropriate $\mathbb{Q}$-form, we may assume that $\mathbf{G}^q$ is a~quasi-split $\mathbb{Q}$-group that
splits over~$L$, and that the $*$-action of $\Gal(L/\mathbb{Q})$ is the same as the $*$-action of
$\Gal(\mathbb{C}/\mathbb{R})$.

Let~$\sigma$ be the nontrivial element of $\Gal(\mathbb{C}/\mathbb{R})$, f\/ix a~representative $a \in
\mathbf{G}^q(\mathbb{C})$ of $\xi(\sigma) \in \overline{\mathbf{G}}{}^q(\mathbb{C})$, and let $z = a {}^\sigma a$.
Since $\sigma^2$ is trivial and $\xi$ is a~$1$-cocycle, we know that~$z$~is trivial in $\overline{\mathbf{G}}{}^q$, so $z
\in Z(\mathbf{G}^q)(\mathbb{C})$.
This implies that~$a$~commutes with~${}^\sigma a$, so~$\sigma$ f\/ixes~$z$, which means $z \in
Z(\mathbf{G}^q)(\mathbb{R})$.

Let $\mathbf{H}$ be the product of the almost simple factors of $\mathbf{G}^q$ that are absolutely almost simple and of
type ${}^{2} A_n$ (more concretely, $\mathbf{H}$ is the product of the factors that are isomorphic to $\SU(k,\ell)$, for
some~$k$ and~$\ell$), let $Z(\mathbf{H})^2 = \{ w^2 \mid w \in Z(\mathbf{H})(\mathbb{C}) \}$, let
$\underline{\mathbf{G}}^q = \mathbf{G}^q/Z(\mathbf{H})^2$, and let $\underline z$ be the image of~$z$ in
$\underline{\mathbf{G}}^q$.
Note that $\underline{\mathbf{G}}^q$ is a~$\mathbb{Q}$-group (since $Z(\mathbf{H})^2$ is a~$\mathbb{Q}$-subgroup
of~$\mathbf{G}^q$).

We claim that we may assume $\underline z \in Z(\underline{\mathbf{G}}^q)(\mathbb{Q})$.
While proving this, we may consider each simple factor individually, so there is no harm in assuming
$\underline{\mathbf{G}}^q$ is almost simple.
This allows us to furthermore assume that $\underline{\mathbf{G}}^q$ is absolutely almost simple.
(Otherwise, since every $\mathbb{C}$-group is split, we could assume $\xi$ is trivial.) Also, since
$|{\Gal(\mathbb{C}/\mathbb{R})}| = 2$, we may assume, by replacing~$a$ with $aw$ for an appropriately chosen $w \in
\langle z \rangle$, that $|\underline z|$~is a~power of~$2$.
Assuming, as we may, that $\underline z$ is nontrivial, this implies that $\mathbf{G}^q$ is not of type ${}^{1,2} E_6$.
Then $Z(\underline{\mathbf{G}}^q)(\mathbb{R}) = Z(\underline{\mathbf{G}}^q)(\mathbb{Q})$.
(If~$\mathbf{G}^q$ is of type ${}^2 A_n$, then the def\/inition of $\underline{\mathbf{G}}^q$ implies
$|Z(\underline{\mathbf{G}}^q)| \le 2$, so every element of $Z(\underline{\mathbf{G}}^q)$ is def\/ined over~$\mathbb{Q}$.
If $\mathbf{G}^q$ is not of this type, then the desired conclusion can be verif\/ied by noting that $Z(\mathbf{G}^q)$ is
either $\boldsymbol\mu_n$, $\boldsymbol\mu_2 \times \boldsymbol\mu_2$, $\Res{L/\mathbb{Q}} \boldsymbol\mu_2$,
$\ResOne{L/\mathbb{Q}} \boldsymbol\mu_2$, or $\ResOne{L/\mathbb{Q}}
\boldsymbol\mu_4$~\cite[p.~332]{PlatonovRapinchukBook}.) This completes the proof of the claim.

The claim of the preceding paragraph implies that the cyclic subgroup $\langle \underline z \rangle$ generated
by~$\underline z$ is def\/ined over~$\mathbb{Q}$.
Hence, the quotient $\widetilde{\mathbf{G}}^q = \mathbf{G}^q/\langle z, Z(\mathbf{H})^2 \rangle$ is a~semisimple
$\mathbb{Q}$-group.
Now, since $a {}^\sigma a = z$ is trivial in $\widetilde{\mathbf{G}}^q$, we know that $\xi$ lifts to a~$1$-cocycle
$\widetilde\xi \colon \Gal(\mathbb{C}/\mathbb{R}) \to \widetilde{\mathbf{G}}^q$.
Then Lemma~\ref{KneserHasse} implies that, after replacing $\widetilde\xi$ with a~cohomologous cocycle, we may assume
$\widetilde\xi$ is the restriction of a~$1$-cocycle $\zeta \colon \Gal(L/\mathbb{Q}) \to \widetilde{\mathbf{G}}^q(L)$.
Let $\mathbf{G}_\mathbb{Q} = {}^\zeta \mathbf{G}^q$, so $\mathbf{G}_\mathbb{Q}$ is a~$\mathbb{Q}$-group that
is~$L$-split, and is isomorphic to~$\mathbf{G}$ over~$\mathbb{R}$.
Also, let $\Lie G$ be the Lie algebra of $\mathbf{G}_\mathbb{Q}$.

To complete the proof, we show that $\Lie G$ is $\mathbb{R}$-universal, by verifying~\eqref{SplitOverQuadIff-weight}.
To this end, let~$\lambda$ be a~$\Gal(\mathbb{C}/\mathbb{R})$-invariant dominant weight, such that $\beta_{\Lie
G,\mathbb{R}}(\lambda) = 0$.
Since $\mathbf{G}^q$~is~$L$-split and the $*$-action of $\Gal(L/\mathbb{Q})$ is the same as the $*$-action of
$\Gal(\mathbb{C}/\mathbb{R})$ (by the choice of the $\mathbb{Q}$-form of~$\mathbf{G}^q$), we know that~$\lambda$~is
invariant under $\Gal(\overline{\mathbb{Q}}/\mathbb{Q})$.

Since~$\zeta$ is a~$1$-cocycle into $\widetilde{\mathbf{G}}^q = \mathbf{G}^q/\langle z, Z(\mathbf{H})^2 \rangle$, we know
that, in the notation of Def\/inition~\ref{BetaDefn} with $F = \mathbb{Q}$, we have $\delta_*[\zeta] \in H^2\bigl(\mathbb{Q};
\langle z, Z(\mathbf{H})^2 \rangle\bigr)$.
Therefore, in order to show that $\beta_{\Lie G, \mathbb{Q}}(\lambda) = \lambda_* \delta_*[\zeta] = 0$, it suf\/f\/ices to
show that~$\lambda$ is trivial on both~$z$ and~$Z(\mathbf{H})^2$.
Note that, in the notation of Def\/inition~\ref{BetaDefn} with $F = \mathbb{R}$, we have $\lambda_* \delta_*[\xi] = \beta_{\Lie G,
\mathbb{R}}(\lambda) = 0$.
Under the natural identif\/ication of $H^2 \bigl(\mathbb{R}; Z(\mathbf{G}^q) \bigr)$ with
\begin{gather*}
\{ w \in Z(\mathbf{G}^q) \mid {}^\sigma w = w \} / \{ w {}^\sigma w \mid w \in Z(\mathbf{G}^q) \},
\end{gather*}
we have $\delta_* [\xi] = [z]$, so this means $\lambda(z) = 1$ (since $\omega \overline{\omega} = 1$ for all $\omega \in
\boldsymbol\mu$).
Furthermore, since the restriction of~$\lambda$ to $Z(\mathbf{G}^q)$ is a~$\Gal(\mathbb{C}/\mathbb{R})$-equivariant
homomorphism, and $Z(\mathbf{H})(\mathbb{R}) = Z(\mathbf{H})$ (cf.~\cite[p.~332]{PlatonovRapinchukBook}), we have
$\lambda \bigl(Z(\mathbf{H}) \bigr) \subseteq \boldsymbol\mu(\mathbb{R}) = \{\pm1\}$, so~$\lambda$ is also trivial on
$Z(\mathbf{H})^2$.
\end{proof}

\begin{proof}[Proof of Corollary~\ref{LattAppl}] Choose an $\mathbb{R}$-universal $\mathbb{Q}$-form $\Lie G_\mathbb{Q}$ of the Lie
algebra~$\Lie G$ of~$G$, let $\overline{\mathbf{G}}$ be the corresponding adjoint $\mathbb{Q}$-group, and let $\Gamma = \{
g \in G \mid \mathop{\mathrm{Ad}} g \in \overline{\mathbf{G}}(\mathbb{Z}) \}$, so~$\Gamma$ is a~lattice in~$G$
(cf.~\cite[Theorem~4.14, p.~220]{PlatonovRapinchukBook}).

Suppose $\rho \colon G \to \GL(n,\mathbb{R})$ is a~f\/inite-dimensional representation of~$G$.
By replacing~$\rho$ with a~conjugate, we may assume $d\rho(\Lie G_\mathbb{Q}) \subseteq \Lie{GL}(n,\mathbb{Q})$ (because
$\Lie G_\mathbb{Q}$ is $\mathbb{R}$-universal).
Then the corresponding representation $\widetilde\rho$ of the universal cover $\mathbf{G}$ of~$\overline{\mathbf{G}}$ is
def\/ined over~$\mathbb{Q}$, so there is a~$\widetilde\rho \bigl(\mathbf{G}(\mathbb{Z}) \bigr)$-invariant
$\mathbb{Z}$-lattice in $\mathbb{Q}^n$~\cite[Remark on p.~173]{PlatonovRapinchukBook}.
Then, since $\widetilde\rho \bigl(\mathbf{G}(\mathbb{Z}) \bigr)$ contains a~f\/inite-index subgroup of $\rho(\Gamma)$,
there is also a~$\rho(\Gamma)$-invariant $\mathbb{Z}$-lattice in $\mathbb{Q}^n$, so $\rho(\Gamma)$ is conjugate to
a~subgroup of $\GL(n,\mathbb{Z})$.
\end{proof}

\section[$\mathbb{R}$-universal absolutely simple Lie algebras]{$\boldsymbol{\mathbb{R}}$-universal absolutely simple Lie algebras}\label{WhichRUnivSect}

This section provides a~classif\/ication of the absolutely simple Lie algebras over~$\mathbb{Q}$ that
are $\mathbb{R}$-universal.
(To say $\Lie G_\mathbb{Q}$ is \emph{absolutely simple} means that $\Lie G_\mathbb{Q} \otimes_\mathbb{Q} \mathbb{C}$ is
simple.
See Section~\ref{NotAbsSect} for a~discussion of $\mathbb{R}$-universal Lie algebras that do not have this property.) We record
a~few observations before proceeding with case-by-case analysis.

\begin{notation}
We assume $\Lie G$ is a~semisimple Lie algebra over a~f\/ield~$F$ of characteristic~$0$, that~$\lambda$~is a~dominant
weight, and the other notation of Def\/inition~\ref{BetaDefn}.
Furthermore, we let $Z^*$ be the f\/inite, abelian group of all homomorphisms from $Z(\mathbf{G}^q)$ to~$\boldsymbol\mu$.
\end{notation}

\begin{prop}[\protect{\cite[Corollary~3.5, \S~4.2, and Lemma~7.4]{Tits-RepLinIrred}}]\label{Beta=D}
Let~$L$ be the center of the algebra $D_{\Lie G, F}(\lambda):= \mathrm{End}_{\Lie G}({}^F \rho_\lambda)$ $($which,
by Schur's lemma, is a~division algebra$)$.
Then:
\renewcommand{\labelenumi}{{\rm (\theenumi)}}
\begin{enumerate}\itemsep=0pt
\item
\label{Beta=D-invt}
$\Gal(\overline{F}/L) = \{ \sigma \in \Gal(\overline{F}/F) \mid \sigma(\lambda) = \lambda \}$ $($for the $*$-action
of $\Gal(\overline{F}/F))$, and
\item $\beta_{\Lie G, L}(\lambda) = [D_{\Lie G, F}(\lambda)]$, after identifying $H^2(L; \boldsymbol\mu)$ with the
Brauer group of~$L$.
\end{enumerate}
\end{prop}

\begin{defn}
We will call $D_{\Lie G, F}(\lambda)$ the \emph{Tits algebra} of~$\Lie G$ corresponding to the weight~$\lambda$
over~$F$.
(However, this name is used in the literature for a~slightly dif\/ferent algebra that is Brauer equivalent to $D_{\Lie G,
F}(\lambda)$~\cite[\S~27A, p.~377]{KnusEtAl-BookOfInvols}.)
\end{defn}

\begin{lem}
\label{RunivIff}
If $\Lie G$ is a~semisimple Lie algebra over~$\mathbb{Q}$, then the following are equivalent:
\renewcommand{\labelenumi}{{\rm \theenumi.}}
\begin{enumerate}\itemsep=0pt
\item
\label{RunivIff-univ}
$\Lie G$ is $\mathbb{R}$-universal.
\item
\label{RunivIff-irred}
${}^\mathbb{Q} \rho_\lambda \otimes_\mathbb{Q} \mathbb{R} \cong {}^\mathbb{R} \rho_\lambda$, for every dominant
weight~$\lambda$.
\item
\label{RunivIff-DivAlg}
$D_{\Lie G, \mathbb{Q}}(\lambda) \otimes_\mathbb{Q} \mathbb{R}$ is a~division algebra, for every dominant
weight~$\lambda$.
\item
\label{RunivIff-D}
$D_{\Lie G, \mathbb{Q}}(\lambda) \otimes_\mathbb{Q} \mathbb{R} \cong D_{\Lie G, \mathbb{R}}(\lambda)$, for every
dominant weight~$\lambda$.
\end{enumerate}
\end{lem}

\begin{proof}
($\ref{RunivIff-univ} \Rightarrow \ref{RunivIff-irred}$) Let~$\rho$ be a~$\mathbb{Q}$-form of ${}^\mathbb{R}
\rho_\lambda$.
Then~$\rho$ must be $\mathbb{Q}$-irreducible (since $\rho \otimes_\mathbb{Q} \mathbb{R} \cong {}^\mathbb{R}
\rho_\lambda$ is $\mathbb{R}$-irreducible) and
\begin{gather*}
\rho \otimes_\mathbb{Q} \mathbb{C} \cong (\rho \otimes_\mathbb{Q} \mathbb{R}) \otimes_\mathbb{R} \mathbb{C} \cong
{}^\mathbb{R} \rho_\lambda \otimes_\mathbb{R} \mathbb{C}
\end{gather*}
has an irreducible summand with highest weight~$\lambda$.
Therefore $\rho \cong {}^\mathbb{Q} \rho_\lambda$.

($\ref{RunivIff-irred} \Rightarrow \ref{RunivIff-univ}$) ${}^\mathbb{Q} \rho_\lambda$ is a~$\mathbb{Q}$-form of
${}^\mathbb{R} \rho_\lambda$.

($\ref{RunivIff-irred} \Leftrightarrow \ref{RunivIff-DivAlg} \Leftrightarrow \ref{RunivIff-D}$) We have
$\mathrm{End}_{\Lie G}({}^\mathbb{Q} \rho_\lambda \otimes_\mathbb{Q} \mathbb{R}) = D_{\Lie G, \mathbb{Q}}(\lambda)
\otimes_\mathbb{Q} \mathbb{R}$.
So ${}^\mathbb{Q} \rho_\lambda \otimes_\mathbb{Q} \mathbb{R}$ is irreducible (and hence equal to~${}^\mathbb{R}
\rho_\lambda$) if and only if $D_{\Lie G, \mathbb{Q}}(\lambda) \otimes_\mathbb{Q} \mathbb{R}$ is a~division algebra (and
hence equal to~$D_{\Lie G, \mathbb{R}}(\lambda)$).
\end{proof}

\begin{rem}
\label{OnlyZStar}
It is immediate from Def\/inition~\ref{BetaDefn} that if the $*$-invariant weights $\lambda_1$ and $\lambda_2$ have the same
restriction to $Z(\mathbf{G}^q)$ (that is, if they represent the same element of~$Z^*$), then $\beta_{\Lie G,
F}(\lambda_1) = \beta_{\Lie G, F}(\lambda_2)$.
Therefore, if~$\lambda$ is a~dominant weight in the root lattice (so~$\lambda$ is trivial on $Z(\mathbf{G}^q)$),
and~$\lambda$ is invariant under the $*$-action of $\Gal(\overline{F}/F)$, then $D_{\Lie G, F}(\lambda) = F$.
\end{rem}

\begin{lem}
\label{SameFixed}
Let  $\Lie G$ be a~semisimple Lie algebra over~$\mathbb{Q}$.
If~$\Lie G$ is $\mathbb{R}$-universal, then $\Lie G$ is inner over some imaginary quadratic extension~$L$
of~$\mathbb{Q}$, and the $*$-action of $\Gal(L/\mathbb{Q})$ is the same as the $*$-action of~$\Gal(\mathbb{C}/\mathbb{R})$.
\end{lem}

\begin{proof}
Assume $\Lie G$ is not inner, for otherwise the desired conclusions are obvious.
Let~$L$ be the (unique) minimal extension of~$\mathbb{Q}$ over which $\Lie G$ becomes inner, and choose some dominant
weight~$\lambda$ that is not f\/ixed by any nontrivial element of~$\Gal(L/\mathbb{Q})$.
We know that~$L$ is the center of~$D_{\Lie G, \mathbb{Q}}(\lambda)$ (see Proposition~\ref{Beta=D}(\ref{Beta=D-invt})) and that $D_{\Lie
G, \mathbb{Q}}(\lambda) \otimes_\mathbb{Q} \mathbb{R}$ is a~division algebra (see Lemma~\ref{RunivIff}(\ref{RunivIff-DivAlg})).
Since~$L$ is not $\mathbb{Q}$, this implies it is an imaginary quadratic extension of~$\mathbb{Q}$.

It now suf\/f\/ices to show that every every $\Gal(\mathbb{C}/\mathbb{R})$-invariant dominant weight~$\lambda$ is also
$\Gal(L/\mathbb{Q})$-invariant.
Suppose not.
By replacing~$\lambda$ with an appropriate positive integer multiple, we may assume it is in the root lattice, so
$D_{\Lie G,\mathbb{R}}(\lambda) = \mathbb{R}$ (see Remark~\ref{OnlyZStar}).
Therefore, we must have $D_{\Lie G,\mathbb{Q}}(\lambda) = \mathbb{Q}$ (see Lemma~\ref{RunivIff}(\ref{RunivIff-D})).
However, this contradicts Proposition~\ref{Beta=D}(\ref{Beta=D-invt}), since~$\lambda$ is not
$\Gal(\overline{\mathbb{Q}}/\mathbb{Q})$-invariant.
\end{proof}

In the remainder of this section, we consider each possible type of absolutely simple Lie algebra over~$\mathbb{Q}$.
The classical types are handled by using the calculations of $D_{\Lie G,F}(\lambda)$ in~\cite[\S~27B,
pp.~378--379]{KnusEtAl-BookOfInvols}, and the answers for ${}^{1,2} E_6$ and $E_7$ follow from observations of Tits
(see Section~\ref{TypeE}).
The remaining types are very easy to deal with:

\begin{lem}
\label{EFGAllUniv}
Every absolutely simple Lie algebra of type $E_8$, $F_4$, or~$G_2$ over~$\mathbb{Q}$ is $\mathbb{R}$-universal.
\end{lem}

\begin{proof}
These types have no outer automorphisms, so the $*$-action must be trivial.
Furthermore, simple groups of these types have trivial center.
Therefore, it is immediate from Remark~\ref{OnlyZStar} that $D_{\Lie G, \mathbb{Q}}(\lambda) = \mathbb{Q}$ for every dominant
weight~$\lambda$.
So Lemma~\ref{RunivIff}(\ref{RunivIff-DivAlg}) implies that $\Lie G$ is $\mathbb{R}$-universal.
\end{proof}

\begin{erratum}
We take this opportunity to correct the statements of Propositions~7.2 and~7.3(a) of~\cite{Morris-RRepsQForms}.
The correct statement of Proposition~7.2 is:
\begin{itemize}\itemsep=0pt
\it
\item[] Suppose $\Lie G_\mathbb{R}$ is a~compact, simple Lie algebra over~$\mathbb{R}$.
There is a $\mathbb{Q}$-form $\Lie G_\mathbb{Q}$ of~$\Lie G_\mathbb{R}$, such that $\Lie G_\mathbb{Q}$ splits
over some quadratic extension of~$\mathbb{Q}$, but is not $\mathbb{R}$-universal, if and only if either
\begin{enumerate}[{\rm (}a{\rm )}]\itemsep=0pt
\item
$\Lie G_\mathbb{R} \cong \Lie{SU}(n)$, for some~$n$ that is divisible by~$4$, or
\item
$\Lie G_\mathbb{R} \cong \Lie{SO}(n)$, for some $n \not\equiv 3,5 \pmod{8}$ $($with $n \ge 6)$.
\end{enumerate}
\end{itemize}
The mistake originates in Proposition~7.3(a) of~\cite{Morris-RRepsQForms}, where~$\ell$~is required to only be odd,
whereas it actually needs to be $\equiv 3 \pmod{4}$.
This means that $\Lie G$, the compact real form of type~$A_\ell$, is isomorphic to $\Lie{SU}(n)$, for some~$n$ that is
divisible by~$4$.
In Proposition~7.2 of~\cite{Morris-RRepsQForms}, it was incorrectly stated that~$n$~only needs to be even, not divisible
by~$4$.
\end{erratum}

\subsection[$\mathbb{R}$-universal Lie algebras of type $A$]{$\boldsymbol{\mathbb{R}}$-universal Lie algebras of type $\boldsymbol{A}$}\label{TypeA}

\begin{prop}
Let $\Lie G = \Lie{SL}_n(D)$, where~$D$ is a~central division algebra over~$\mathbb{Q}$.
Then $\Lie G$ is $\mathbb{R}$-universal if and only if~$D$ is either~$\mathbb{Q}$ or a~quaternion algebra that
does not split over~$\mathbb{R}$.
\end{prop}

\begin{proof}\looseness=1
Let~$d$ be the degree of~$D$ over~$\mathbb{Q}$, and let~$\lambda$ be the highest weight of the standard representation
of $\Lie{SL}_{dn}(\mathbb{C})$, so $\{\lambda^i\}_{i=1}^{dn}$ is a~set of representatives of~$Z^*$.
Then $D_{\Lie G,\mathbb{Q}}(\lambda^i)$ and $D_{\Lie G,\mathbb{R}}(\lambda^i)$ are Brauer equivalent to the~$i$-fold
tensor products $D^{\otimes i}$ and $D_{\Lie G, \mathbb{R}}(\lambda)^{\otimes i}$,
respectively~\cite[p.~378]{KnusEtAl-BookOfInvols}.

If $\Lie G$ is $\mathbb{R}$-universal, then, by taking $i = 1$ and noting that $D_{\Lie G,\mathbb{R}}(\lambda)$ is
either~$\mathbb{R}$ or the quaternion algebra~$\mathbb{H}$, we see from Lemma~\ref{RunivIff}(\ref{RunivIff-D}) that~$D$ must be
either $\mathbb{Q}$ or a~quaternion algebra that does not split over~$\mathbb{R}$.

Conversely, suppose~$D$ is either $\mathbb{Q}$ or a~quaternion algebra that does not split over~$\mathbb{R}$.
In either case, $D^2$ is Brauer equivalent to~$\mathbb{Q}$, so $D_{\Lie G,\mathbb{Q}}(\lambda^i) \otimes_\mathbb{Q}
\mathbb{R}$ is $\mathbb{Q} \otimes_\mathbb{Q} \mathbb{R} = \mathbb{R}$ if~$i$ is even, and it is $D \otimes_\mathbb{Q}
\mathbb{R} = \mathbb{H}$ if~$i$ is odd.
This is a~division algebra for every~$i$, so we see from Lemma~\ref{RunivIff}(\ref{RunivIff-DivAlg}) that $\Lie G$ is
$\mathbb{R}$-universal.
\end{proof}

\begin{prop}
Let $\Lie G = \Lie{SU}_n(B; D, \tau)$, where
\begin{itemize}\itemsep=0pt
\item $D$ is a~division algebra that is central over a~quadratic extension~$L$ of~$\mathbb{Q}$,
\item $\tau$ is an anti-involution of~$D$ that is nontrivial on~$L$, and
\item $B$ is an invertible~$\tau$-Hermitian matrix in $\Mat_{n\times n}(D)$.
\end{itemize}
Then $\Lie G$ is $\mathbb{R}$-universal if and only if $D = L$ is an imaginary quadratic extension
of~$\mathbb{Q}$, and either~$n$~is odd or $(-1)^{n/2} \det B$ is either negative or the norm of some element of~$L$.
\end{prop}

\begin{proof}
We prove only ($\Rightarrow$), but the argument is reversible.
Since $\Lie G$ is outer, Lemma~\ref{SameFixed} implies that $\Lie G \otimes_\mathbb{Q} \mathbb{R}$ is also outer, and the
quadratic extension~$L$ is imaginary.

Let~$d$ be the degree of~$D$ over~$L$, and let~$\lambda$ be the highest weight of the standard representation of
$\Lie{SL}_{dn}(\mathbb{C})$, so $\{\lambda^i\}_{i=1}^{dn}$ is a~set of representatives of~$Z^*$.

The Tits algebra of the natural representation $\rho \colon \Lie G \hookrightarrow \Mat_{n \times n}(D)$
is~$D$~\cite[p.~378]{KnusEtAl-BookOfInvols}.
Since~$D$ splits over~$\mathbb{R}$ (recall that the center~$L$ is an imaginary extension),
we have $D_{\Lie G, \mathbb{R}}(\lambda) = \mathbb{C}$, so we conclude from Lemma~\ref{RunivIff}(\ref{RunivIff-D})
that $D = L$.

If~$\lambda$ is any weight that is not f\/ixed by the $*$-action, then $D_{\Lie G, \mathbb{Q}}(\lambda) = L$ and $D_{\Lie
G, \mathbb{R}}(\lambda) = \mathbb{C}$~\cite[p.~378]{KnusEtAl-BookOfInvols}, so $D_{\Lie G, \mathbb{Q}}(\lambda)
\otimes_\mathbb{Q} \mathbb{R} = D_{\Lie G, \mathbb{R}}(\lambda)$, as specif\/ied in Lemma~\ref{RunivIff}(\ref{RunivIff-D}).
Hence, such weights (which are all of the weights when~$n$ is odd) do not place any further restriction on~$\Lie G$.

Suppose~$n$ is even.
Any weight that is f\/ixed by the $*$-action (and is not in the root lattice) is congruent to $\lambda^{n/2}$, modulo the
root lattice.
Write $L = \mathbb{Q} \bigl[\sqrt{a} \bigr]$ and let $b = (-1)^{n/2} \det B$.
Then, for~$F$ either~$\mathbb{Q}$ or~$\mathbb{R}$, $D_{\Lie G, F}(\lambda^{n/2})$ is Brauer equivalent to the quaternion
algebra $(a,b)_F$~\cite[p.~378 and Corollary~10.35 on p.~131)]{KnusEtAl-BookOfInvols}.
This is trivial in the Brauer group if and only if it is split, which means that~$b$ is the norm of some element of $F
\bigl[\sqrt{a} \bigr]$.

Since $\Lie G$ is $\mathbb{R}$-universal, we see from Lemma~\ref{RunivIff}(\ref{RunivIff-D}) that either $D_{\Lie G,
\mathbb{R}}(\lambda^{n/2}) = \mathbb{H}$ (which, by the preceding paragraph, means that~$b$ is not a~norm in $\mathbb{R}
\bigl[\sqrt{a} \bigr] = \mathbb{C}$, so $b < 0$) or $D_{\Lie G, \mathbb{Q}}(\lambda^{n/2}) = \mathbb{Q}$ (which means
that~$b$~is a~norm in $\mathbb{Q} \bigl[\sqrt{a} \bigr] = L$).
\end{proof}

\subsection[$\mathbb{R}$-universal Lie algebras of type $B$]{$\boldsymbol{\mathbb{R}}$-universal Lie algebras of type $\boldsymbol{B}$}\label{TypeB}

\begin{notation}
Let $F \in \{\mathbb{Q}, \mathbb{R}\}$.
Any symmetric matrix $B \in \GL_k(F)$ determines a~nondegenerate quadratic form on~$F^k$.
We use $\Cliff^0_F(B)$ to denote the corresponding even Clif\/ford algebra~\cite[p.~104]{Lam-IntroQuadForms}.
\end{notation}

It is well known that $\Cliff^0_F(B)$ is either a~simple algebra or the direct sum of two isomorphic simple algebras
over~$F$~\cite[Theorems~2.4 and~2.5, p.~110]{Lam-IntroQuadForms}.
If~$B$~has been diagonalized, then it is straightforward to determine whether this simple algebra is split (in which
case, we also say that~$\Cliff^0_F(B)$ is split).
Namely, the simple algebra is Brauer equivalent to a~quaternion algebra that can be calculated from the eigenvalues
of~$B$ (cf.~\cite[Corollary~3.14, p.~117]{Lam-IntroQuadForms}).

\begin{eg}[\protect{\cite[pp.~122--125, and Corollary~2.10, p.~112]{Lam-IntroQuadForms}}] Assume~$B$ is a~symmetric matrix in
$\GL_k(\mathbb{R})$ with exactly~$p$ positive eigenvalues (including multiplicity, so $k - p$ is the number of negative
eigenvalues).
Then $\Cliff^0_\mathbb{R}(B)$ is not split if and only if $2p - k$ is congruent to~$3$,~$4$, or~$5$, modulo~$8$.
\end{eg}

\begin{prop}
Let $\Lie G = \Lie{SO}_n(B; \mathbb{Q})$, where~$B$~is a~symmetric matrix in $\GL_n(\mathbb{Q})$, and~$n$~is odd.
Then $\Lie G$ is $\mathbb{R}$-universal if and only if either $\Cliff^0_\mathbb{Q}(B)$ is split or
$\Cliff^0_\mathbb{R}(B)$ is not split.
\end{prop}

\begin{proof}
Since the center of $\mathbf{Spin}_n$ has order~$2$ when~$n$~is odd, there is only one Tits algebra to consider
(see Remark~\ref{OnlyZStar}), and it is Brauer equivalent to $\Cliff^0_F(B)$~\cite[p.~378]{KnusEtAl-BookOfInvols}.
\end{proof}

\subsection[$\mathbb{R}$-universal Lie algebras of type $C$]{$\boldsymbol{\mathbb{R}}$-universal Lie algebras of type $\boldsymbol{C}$}\label{TypeC}

\begin{eg}
For any~$n$, the Lie algebra $\Lie{Sp}_{2n}(\mathbb{Q})$ is $\mathbb{R}$-universal.
(More generally, every $\mathbb{Q}$-split semisimple Lie algebra is $\mathbb{R}$-universal, because it is clear from
Def\/inition~\ref{BetaDefn} that $\beta_{\Lie G, F}(\lambda)$ is always trivial for~$F$-split Lie algebras.)
\end{eg}

\begin{prop}
Let $\Lie G = \Lie{SU}_n(B; D, \tau)$, where
\begin{itemize}\itemsep=0pt
\item $D$ is a~quaternion division algebra over~$\mathbb{Q}$,
\item $\tau$ is the reversion anti-involution of~$D$, and
\item $B$ is an invertible~$\tau$-Hermitian matrix in $\Mat_{n\times n}(D)$.
\end{itemize}
Then $\Lie G$ is $\mathbb{R}$-universal if and only if~$D$~does not split over~$\mathbb{R}$.
\end{prop}

\begin{proof}
Since the center of $\mathbf{Sp}_n$ has order~$2$, there is only one Tits algebra to consider (see Remark~\ref{OnlyZStar}).
This Tits algebra is~$D$ (over~$\mathbb{Q}$)~\cite[p.~378]{KnusEtAl-BookOfInvols}, so the desired conclusion is
immediate from Lemma~\ref{RunivIff}(\ref{RunivIff-DivAlg}).
\end{proof}

\subsection[$\mathbb{R}$-universal Lie algebras of type $D$]{$\boldsymbol{\mathbb{R}}$-universal Lie algebras of type $\boldsymbol{D}$}\label{TypeD}

Lemma~\ref{SameFixed} implies that triality forms are not $\mathbb{R}$-universal.
Therefore, all absolutely simple $\mathbb{R}$-universal Lie algebras of type~$D$ are described in either Proposition~\ref{so2k} or
Lemma~\ref{SO(quat)}.

\begin{prop}
\label{so2k}
Let $\Lie G = \Lie{SO}_{2k}(B; \mathbb{Q})$, with $k \ge 3$, for some symmetric $B \in \GL_{2k}(\mathbb{Q})$.
Then~$\Lie G$ is $\mathbb{R}$-universal if and only if either
\renewcommand{\labelenumi}{{\rm (\theenumi)}}
\begin{enumerate}\itemsep=0pt
\item $\Cliff^0_\mathbb{Q}(B)$ is split, and $(-1)^k \det B$ is either negative or a~square in~$\mathbb{Q}$, or
\item $\Cliff^0_\mathbb{Q}(B)$ does not split over~$\mathbb{R}$, and $(-1)^k \det B$ is a~square in~$\mathbb{Q}$.
\end{enumerate}
\end{prop}

\begin{proof}
As in~\cite[\S~27B, type~$D_n$, p.~379]{KnusEtAl-BookOfInvols}, let $\lambda$, $\lambda_+$, $\lambda_-$ be dominant weights
that represent the three nonzero elements of~$Z^*$, with~$\lambda$ being the highest weight of the natural
representation of~$\Lie G$ on $\mathbb{C}^{2k}$.
Then $D_{\Lie G, \mathbb{Q}}(\lambda) = \mathbb{Q}$ is trivial~\cite[p.~379]{KnusEtAl-BookOfInvols}.

Suppose, f\/irst, that $\Lie G$ is an inner form, which means that $(-1)^k \det B$ is a~square in~$\mathbb{Q}$.
Then $\Cliff^0_\mathbb{Q}(B)$ is a~direct sum of two algebras $C^+$ and~$C^-$ that are Brauer equivalent to the full
Clif\/ford algebra~\cite[Theorem~2.5(3), p.~110]{Lam-IntroQuadForms} (and are therefore Brauer equivalent to a~quaternion
algebra).
Furthermore, $D_{\Lie G, \mathbb{Q}}(\lambda_{\pm})$ is Brauer equivalent
to~$C^\pm$~\cite[p.~379]{KnusEtAl-BookOfInvols}.
Therefore, Lemma~\ref{RunivIff}(\ref{RunivIff-DivAlg}) shows that an inner form $\Lie G$ is:
\begin{itemize}\itemsep=0pt
\item automatically $\mathbb{R}$-universal, when $\Cliff^0_\mathbb{Q}(B)$ does not split over~$\mathbb{R}$, but
\item $\mathbb{R}$-universal if and only if $\Cliff^0_\mathbb{Q}(B)$ is split, when $\Cliff^0_\mathbb{Q}(B)$ splits
over~$\mathbb{R}$.
\end{itemize}

Assume, now, that $\Lie G$ is an outer form.
In order for $\Lie G$ to be $\mathbb{R}$-universal, $\Lie G$ must remain outer over~$\mathbb{R}$ (see Lemma~\ref{SameFixed}),
which means $(-1)^k \det B < 0$.
Let~$L$ be the (unique, imaginary) quadratic extension of~$\mathbb{Q}$ over which $\Lie G$ becomes inner.
Then $D_{\Lie G, \mathbb{Q}}(\lambda_{\pm})$ is Brauer equivalent to
$\Cliff^0_\mathbb{Q}(B)$~\cite[p.~379]{KnusEtAl-BookOfInvols}, which is central simple over~$L$.
Since~$B$ splits over~$\mathbb{R}$ (recall that the quadratic extension~$L$ is
imaginary), Lemma~\ref{RunivIff}(\ref{RunivIff-DivAlg}) implies that $\Lie G$ is $\mathbb{R}$-universal if and only if
$\Cliff^0_\mathbb{Q}(B)$ is split.
\end{proof}

\begin{notation}
The notion of even Clif\/ford algebra was extended to the situation of Lem\-ma~\ref{SO(quat)} below by N.~Jacobson.
(A construction can be found in~\cite[\S~8B, pp.~91ff]{KnusEtAl-BookOfInvols}.) We will denote this algebra by $C^0_D(B,
\tau_r)$.
\end{notation}

\begin{lem}
\label{SO(quat)}
Let $\Lie G = \Lie{SU}_k(B; D, \tau_r)$, where
\begin{itemize}\itemsep=0pt
\item $D$ is a~quaternion division algebra over~$\mathbb{Q}$,
\item $\tau_r$ is the reversion anti-involution, and
\item $B$~is a~$\tau_r$-Hermitian matrix in $\GL_k(D)$.
\end{itemize}
Then $\Lie G$ is $\mathbb{R}$-universal if and only if~$D$ does not split over~$\mathbb{R}$, and either
\renewcommand{\labelenumi}{{\rm (\theenumi)}}
\begin{enumerate}\itemsep=0pt
\item $k$ is even and the reduced norm of~$B$ $($calculated in the algebra $\Mat_{k \times k}(D))$ is
a~square in~$\mathbb{Q}$, or
\item $k$ is odd and $C^0_D(B, \tau_r)$ is split.
\end{enumerate}
\end{lem}

\begin{proof}
As in the proof of Proposition~\ref{so2k}, let $\lambda$, $\lambda_+$, $\lambda_-$ be dominant weights that represent the three nonzero
elements of~$Z^*$, with~$\lambda$ being the highest weight of the natural representation of~$\Lie G$ on
$(D\otimes_\mathbb{Q} \mathbb{C})^k$.

($\Rightarrow$) We have $D_{\Lie G, \mathbb{Q}}(\lambda) = D$~\cite[p.~379]{KnusEtAl-BookOfInvols}, so~$D$ does not
split over~$\mathbb{R}$ (see Lemma~\ref{RunivIff}(\ref{RunivIff-DivAlg})).
(Note that this implies $\Lie G \otimes_\mathbb{Q} \mathbb{R} \cong \Lie{SO}_k(\mathbb{H})$.)

Suppose~$k$ is even, so $\Lie G \otimes_\mathbb{Q} \mathbb{R}$ is inner.
Then $\Lie G$ must also be inner (see Lemma~\ref{SameFixed}).
Since~$k$ is even, this means that the reduced norm of~$B$ is a~square in~$\mathbb{Q}$.

Suppose~$k$ is odd.
Then $\Lie G \otimes_\mathbb{Q} \mathbb{R}$ is outer, so the weights $\lambda_+$ and~$\lambda_-$ are not f\/ixed by the
$*$-action.
Therefore $D_{\Lie G, \mathbb{R}}(\lambda_{\pm}) = \mathbb{C}$ (see Proposition~\ref{Beta=D}).
Since $D_{\Lie G, \mathbb{Q}}(\lambda_{\pm})$ is Brauer equivalent to $C^0_D(B,
\tau_r)$~\cite[p.~379]{KnusEtAl-BookOfInvols}, we conclude that $\Lie G$ is $\mathbb{R}$-universal if and only if this
algebra is split (see Lemma~\ref{RunivIff}(\ref{RunivIff-D})).

($\Leftarrow$) The proof for odd~$k$ is reversible, so let us assume~$k$ is even.
Then $D_{\Lie G, \mathbb{R}}(\lambda_{\pm}) = \mathbb{H}$ is nontrivial.
(For example, $D_{\Lie G, \mathbb{R}}(\lambda_{\pm})$ can be calculated by using~\cite[\S~5.5]{Tits-RepLinIrred}.)
Therefore $D_{\Lie G, \mathbb{Q}}(\lambda_{\pm})$ does not split over~$\mathbb{R}$.
Since $D_{\Lie G, \mathbb{Q}}(\lambda_{\pm})$ is a~quaternion algebra (because the center of $\mathbf{Spin}_{2k}$ has
exponent~$2$ when~$k$~is even), we conclude that $D_{\Lie G, \mathbb{Q}}(\lambda_{\pm}) \otimes_\mathbb{Q} \mathbb{R}
\cong D_{\Lie G, \mathbb{R}}(\lambda_{\pm})$, as desired.
\end{proof}

\subsection[$\mathbb{R}$-universal Lie algebras of type $E_6$ and $E_7$]{$\boldsymbol{\mathbb{R}}$-universal Lie algebras of type $\boldsymbol{E_6}$ and $\boldsymbol{E_7}$}\label{TypeE}

It was pointed out in Lemma~\ref{EFGAllUniv} that every absolutely simple Lie algebra of type $E_8$, $F_4$, or~$G_2$ is
$\mathbb{R}$-universal.
In addition, the classical types were discussed in Sections~\ref{TypeA}--\ref{TypeD}.
Therefore, the only types that remain are $E_6$ and~$E_7$.

\begin{defn}[\protect{\cite[p.~649]{Tits-StronglyInner}}] In the notation of Def\/inition~\ref{BetaDefn}, a~Lie algebra $\Lie G$ over a~f\/ield~$F$ is
\emph{strongly inner} if $\mathbf{G}^q$ is split and the cohomology class $[\xi] \in H^1(F; \overline{\mathbf{G}}{}^q)$ is
the image of a~cohomology class in $H^1(F; \mathbf{G}{}^q)$ (under the map induced by the natural homomorphism
$\mathbf{G}{}^q \to \overline{\mathbf{G}}{}^q$).
This condition on $[\xi]$ is equivalent to requiring that $\delta_* [\xi] = 0$ in $H^2\bigl(F; Z(\mathbf{G}^q) \bigr)$,
or, in other words, that $D_{\Lie G, F}(\lambda) = F$ for every dominant weight~$\lambda$.
\end{defn}

\begin{prop}
Let $\Lie G$ be an absolutely simple Lie algebra over~$\mathbb{Q}$ of type~$E_6$.
Then $\Lie G$ is $\mathbb{R}$-universal if and only if either
\renewcommand{\labelenumi}{{\rm (\theenumi)}}
\begin{enumerate}\itemsep=0pt
\item $\Lie G$ is strongly inner $($that is, of type ${}^1 E^{28}_{6,2}$ or ${}^1 E^{0}_{6,6})$, or
\item $\Lie G$ is an outer form that splits over an imaginary quadratic extension of~$\mathbb{Q}$.
\end{enumerate}
\end{prop}

\begin{proof}
Assume, f\/irst, that $\Lie G$ is inner.
Since the center of any simply connected, almost simple group of type~$E_6$ is cyclic of prime order (namely, it is of
order~$3$), we see that $\Lie G$ is $\mathbb{R}$-universal if and only if either $\Lie G$ is strongly inner
(over~$\mathbb{Q}$) or $\Lie G \otimes_\mathbb{Q} \mathbb{R}$ is not strongly inner.
From~\cite[\S~6.4.5]{Tits-RepLinIrred}, we see that Lie algebras of type ${}^1 E^{28}_{6,2}$ or ${}^1 E^{0}_{6,6}$ are
always strongly inner, and those of type ${}^1 E^{16}_{6,2}$ are never strongly inner.
In addition, a~Lie algebra of type ${}^1 E^{78}_{6,0}$ over~$\mathbb{Q}$ cannot be strongly inner~\cite[Propositions~4
and~5]{Tits-StronglyInner}.
Since $\Lie G \otimes_\mathbb{Q} \mathbb{R}$ must be of type ${}^1 E^{28}_{6,2}$ or ${}^1 E^{0}_{6,6}$, it is strongly
inner.
Hence, $\Lie G$ is $\mathbb{R}$-universal if and only if it is strongly inner.

Suppose, now, that $\Lie G$ is outer, and let~$L$ be the unique quadratic extension of~$\mathbb{Q}$ over which $\Lie
G$ is inner.
From Lemma~\ref{SameFixed}, we know that if $\Lie G$ is $\mathbb{R}$-universal, then it must remain outer over~$\mathbb{R}$,
so~$L$~is an imaginary extension of~$\mathbb{Q}$.
Let~$\lambda$ be a~weight that represents a~nontrivial element of~$Z^*$.
Then~$\lambda$ is not f\/ixed by the $*$-action, so $D_{\Lie G, \mathbb{R}}(\lambda) = D_{\Lie G, \mathbb{C}}(\lambda) =
\mathbb{C}$ and $D_{\Lie G, \mathbb{Q}}(\lambda) = D_{\Lie G, L}(\lambda)$ (see Proposition~\ref{Beta=D}).
So $\Lie G$ is $\mathbb{R}$-universal if and only if $\Lie G \otimes_\mathbb{Q} L$ is strongly inner.

Since~$L$ is an imaginary extension, we know that $\Lie G$ splits at the inf\/inite place of~$L$.
Then, by inspection of the possible Tits indices of type ${}^1 E_6$ over a~nonarchimedean local
f\/ield~\cite[p.~58]{Tits-Classification}, we see that the central vertex of the Tits index is circled at every place, so
it must be circled over~$L$~\cite[Satz~4.3.3]{Harder-BerichtNeuere}.
Therefore, $\Lie G$ must be either split or of type ${}^1 E_{6,2}^{16}$ over~$L$.
From~\cite[\S~6.4.5]{Tits-RepLinIrred}, we see that the form of type ${}^1 E_{6,2}^{16}$ is not strongly inner.
Therefore, $\Lie G$ is $\mathbb{R}$-universal if and only if $\Lie G \otimes_\mathbb{Q} L$ is split.
\end{proof}

\begin{prop}
Let $\Lie G$ be an absolutely simple Lie algebra over~$\mathbb{Q}$ of type~$E_7$.
Then $\Lie G$ is $\mathbb{R}$-universal if and only if either
\renewcommand{\labelenumi}{{\rm (\theenumi)}}
\begin{enumerate}\itemsep=0pt
\item $\Lie G$ is of type $E^{28}_{7,3}$ or $E^0_{7,7}$ $($over~$\mathbb{Q})$, or
\item $\Lie G \otimes_\mathbb{Q} \mathbb{R}$ is of type $E^{133}_{7,0}$, $E^{31}_{7,2}$, or $E^9_{7,4}$.
\end{enumerate}
\end{prop}

\begin{proof}
Since a~simply connected group of type~$E_7$ has a~center of order~$2$, there is only one Tits algebra to consider.
For $\Lie G$ to be $\mathbb{R}$-universal, this algebra needs to either be trivial over~$\mathbb{Q}$ or nontrivial
over~$\mathbb{R}$.
Tits~\cite[\S~6.5.5]{Tits-RepLinIrred} points out that it is trivial (over any f\/ield) for types $E^{28}_{7,3}$ and
$E^0_{7,7}$, but nontrivial (over any f\/ield) for types $E^{31}_{7,2}$ and $E^9_{7,4}$.
Since there are no strongly inner anisotropic groups of type~$E_7$ over~$\mathbb{Q}$ or~$\mathbb{R}$~\cite[Propositions~4
and~5]{Tits-StronglyInner}, the Tits algebra is also nontrivial for the anisotropic Lie algebra $E^{133}_{7,0}$ (over
the f\/ields of interest to us).
\end{proof}

\section[$\mathbb{R}$-universal Lie algebras that are not absolutely simple]{$\boldsymbol{\mathbb{R}}$-universal Lie algebras that are not absolutely simple}\label{NotAbsSect}

Section~\ref{WhichRUnivSect} lists the $\mathbb{R}$-universal Lie algebras that are absolutely simple.
It is easy to describe the rest of the simple ones:

\begin{prop}
Let $\Lie G$ be a~simple Lie algebra over~$\mathbb{Q}$ that is not absolutely simple.
Then $\Lie G$ is $\mathbb{R}$-universal if and only if $\Lie G = \Res{L/\mathbb{Q}} \Lie G'$ is the restriction of
scalars of a~strongly inner, absolutely simple Lie algebra~$\Lie G'$ over an imaginary quadratic extension~$L$
of~$\mathbb{Q}$.
\end{prop}

\begin{proof}
($\Rightarrow$) We have $\Lie G = \Res{L/\mathbb{Q}} \Lie G'$, for some absolutely simple Lie algebra~$\Lie G'$ over
some f\/inite extension~$L$ of~$\mathbb{Q}$~\cite[\S~3.1.2]{Tits-Classification}.
From Lemma~\ref{SameFixed}, we see that if $\Lie G$ is $\mathbb{R}$-universal, then the extension~$L$ must be imaginary
quadratic.

Let $\lambda'$ be a~nonzero dominant weight of~$\Lie G'$, and let $\lambda = \Res{L/\mathbb{Q}} \lambda'$ be the
corresponding weight of~$\Lie G$.
Then~$\lambda$ is not f\/ixed by the $*$-action, so we have $D_{\Lie g, F}(\lambda) = D_{\Lie g, L}(\lambda) = D_{\Lie g',
L}(\lambda)$.
(The f\/inal equality is because $\Lie G = \Res{L/\mathbb{Q}} \Lie G'$.) On the other hand, since all Lie algebras are
split over~$\mathbb{C}$, we have $D_{\Lie g, \mathbb{R}}(\lambda) = D_{\Lie g, \mathbb{C}}(\lambda) = \mathbb{C}$.
Hence, we see from Lemma~\ref{RunivIff}(\ref{RunivIff-D}) that if $\Lie G$ is $\mathbb{R}$-universal, then $D_{\Lie g', L}(\lambda)
= L$.
Since this is true for every dominant weight~$\lambda$ of~$\Lie G'$, we conclude that $\Lie G'$ is strongly inner.

($\Leftarrow$) Note that, since $\Lie G'$ is strongly inner, the simply connected group~$\mathbf{G}'$ corresponding
to~$\Lie G'$ may be written as $\mathbf{G}' = {}^{\zeta} \mathbf{G}^s$, where $\mathbf{G}^s$ is split (and simply
connected) and $\zeta \in H^1(L; \mathbf{G}^s)$.
By Restriction of Scalars, then $\mathbf{G} = {}^{\xi} \mathbf{G}^q$, where $\mathbf{G}^q = \Res{L/\mathbb{Q}}
\mathbf{G}^s$ is quasi-split and $\xi \in H^1(\mathbb{Q}; \mathbf{G}^q)$.
Then, in the notation of Def\/inition~\ref{BetaDefn}, we have $\delta_*[\xi] = 0$, so $\beta_{\Lie G, F}(\lambda) = 0$ for every
$*$-invariant weight~$\lambda$.
It is therefore easy to see that $\Lie G$ is $\mathbb{R}$-universal.
\end{proof}

We also brief\/ly describe the considerations involved in constructing semisimple $\mathbb{R}$-universal Lie algebras from
simple ones:

\begin{prop}
A~direct sum $\Lie G = \Lie G_1 \oplus \dots \oplus \Lie G_r$ of simple Lie algebras over~$\mathbb{Q}$ is
$\mathbb{R}$-universal if and only if:
\renewcommand{\labelenumi}{{\rm (\theenumi)}}
\begin{enumerate}\itemsep=0pt
\item
\label{SS-Gi}
each $\Lie G_i$ is $\mathbb{R}$-universal,
\item
\label{SS-star}
either $\Lie G$ is inner, or it becomes inner over an imaginary quadratic extension~$L$ of~$\mathbb{Q}$, and the
$*$-action of $\Gal(\mathbb{C}/\mathbb{R})$ is the same as the $*$-action of $\Gal(L/\mathbb{Q})$.
\item
\label{SS-Di=Dj}
we have $D_{\Lie G_i, \mathbb{Q}}(\lambda_i) \cong D_{\Lie G_j, \mathbb{Q}}(\lambda_j)$ whenever $\lambda_i$
and~$\lambda_j$ are $*$-invariant dominant weights of $G_i$ and~$G_j$, such that $D_{\Lie G_i, \mathbb{Q}}(\lambda_i)
\neq \mathbb{Q}$ and $D_{\Lie G_j, \mathbb{Q}}(\lambda_j) \neq \mathbb{Q}$, and
\item
\label{SS-Dsplit}
either $\Lie G$ is inner, or it becomes inner over an imaginary quadratic extension~$L$ of~$\mathbb{Q}$, such that
$D_{\Lie G_i, \mathbb{Q}}(\lambda_i)$ splits over~$L$, for every dominant weight~$\lambda_i$ of every $\Lie G_i$ that is
inner $($over~$\mathbb{Q})$.
\end{enumerate}
\end{prop}

\begin{proof}
We prove only ($\Rightarrow$), but the argument is reversible.
(We use $\sim$ to denote ``is Brauer equivalent to''.)

(\ref{SS-Gi}) Any representation $\rho_i$ of~$\Lie G_i$ extends to a~representation of $\Lie G$ (that is $0$ on the
other simple factors), so it is clear that $\Lie G_i$ must be $\mathbb{R}$-universal.

(\ref{SS-star}) See Lemma~\ref{SameFixed}.
(This means that all $\Lie G_i$ become inner over the same quadratic extension~$L$.)

(\ref{SS-Di=Dj}) We have
\begin{gather*}
D_{\Lie G, \mathbb{R}}(\lambda_i + \lambda_j) \sim D_{\Lie G, \mathbb{R}}(\lambda_i) \otimes_\mathbb{R} D_{\Lie G,
\mathbb{R}}(\lambda_j) \\
\hphantom{D_{\Lie G, \mathbb{R}}(\lambda_i + \lambda_j)}{}
\cong \bigl(D_{\Lie G, \mathbb{Q}}(\lambda_i) \otimes_\mathbb{Q} \mathbb{R} \bigr)
\otimes_\mathbb{R} \bigl(D_{\Lie G, \mathbb{Q}}(\lambda_i) \otimes_\mathbb{Q} \mathbb{R} \bigr) \cong \mathbb{H}
\otimes_\mathbb{R} \mathbb{H} \sim \mathbb{R}.
\end{gather*}
Since $\Lie G$ is $\mathbb{R}$-universal, this implies $D_{\Lie G, \mathbb{Q}}(\lambda_i + \lambda_j) \sim\mathbb{Q}$,
from which we conclude that $D_{\Lie G, \mathbb{Q}}(\lambda_i) \otimes_\mathbb{Q} D_{\Lie G, \mathbb{Q}}(\lambda_j) \sim
\mathbb{Q}$, so $D_{\Lie G, \mathbb{Q}}(\lambda_i) \cong D_{\Lie G, \mathbb{Q}}(\lambda_j)$ (since these are quaternion
algebras).

(\ref{SS-Dsplit}) Let $\lambda_j$ be a~dominant weight of some (outer) simple factor $G_j$ that is not f\/ixed by the
$*$-action.
Then $\lambda_i + \lambda_j$ is not f\/ixed by the $*$-action, so
\begin{gather*}
D_{\Lie G, \mathbb{Q}}(\lambda_i + \lambda_j) = D_{\Lie G, L}(\lambda_i + \lambda_j) \sim D_{\Lie G_i, L}(\lambda_i)
\otimes_L D_{\Lie G_j, L}(\lambda_j) = D_{\Lie G_i, L}(\lambda_i) \otimes_L L \cong D_{\Lie G_i, L}(\lambda_i).
\end{gather*}
Since $D_{\Lie G, \mathbb{R}}(\lambda_i + \lambda_j) = D_{\Lie G, \mathbb{C}}(\lambda_i + \lambda_j) = \mathbb{C}$, and
$\Lie G$ is $\mathbb{R}$-universal, we conclude from Lemma~\ref{RunivIff}(\ref{RunivIff-D}) that $D_{\Lie G_i, L}(\lambda_i) = L$.
Since $D_{\Lie G, \mathbb{Q}}(\lambda_i) \otimes_\mathbb{Q} L \sim D_{\Lie G, L}(\lambda_i)$, this means that $D_{\Lie
G, \mathbb{Q}}(\lambda_i)$ splits over~$L$.
\end{proof}

\subsection*{Acknowledgements}
It is a~pleasure to thank V.~Chernousov for a~very helpful discussion about Tits algebras of special orthogonal groups,
A.~Rapinchuk for explaining how to prove~Lemma~\ref{KneserHasse}, and the anonymous referees for numerous very insightful
comments on a~previous version of this manuscript, inclu\-ding some important corrections.

\pdfbookmark[1]{References}{ref}
\LastPageEnding

\end{document}